\title{Analysis of disk scheduling, increasing subsequences and space-time geometry}
\author{Eitan Bachmat 
\thanks{
Department of computer science,
Ben-Gurion University, Email:ebachmat@cs.bgu.ac.il.
 Research partially supported by an IBM faculty award and a research grant from Seagate Corp.}}
\date{}
\documentstyle [11pt]{article}
\newtheorem{theo}{Theorem}[section]

\newtheorem{coro}[theo]{Corollary}

\begin{document}
\maketitle

\begin{abstract}

We consider the problem of estimating the average tour length of the 
asymmetric TSP arising from the 
disk scheduling problem with a linear seek function  
and a probability distribution on the location of I/O requests.  
The optimal disk scheduling algorithm of Andrews, Bender and Zhang is interpreted as a simple peeling process on points in a 2 dimensional space-time
w.r.t the causal structure. The  
patience sorting algorithm for finding the longest increasing subsequence in a permutation can be 
given a
similar interpretation. Using this interpretation 
we show that the optimal tour length is
the length of the maximal curve with respect to a Lorentzian metric 
on the surface of the disk drive. This length can be computed explicitly
in some interesting cases. When the probability distribution is assumed uniform we provide 
finer asymptotics for the tour length.
The interpretation also provides a  better understanding of patience sorting and allows us to extend a result of Aldous and Diaconis on pile sizes 

\end{abstract}

\section{introduction}
Modern disk drives have the ability to queue incoming read and write  
requests and to service them in an out of order fashion. 
In the {\it Batched disk scheduling} problem we are given a batch of 
$n$ queued requests and we wish to service them in an order which minimizes 
the total service time, or equivalently, in an order which minimizes the number of disk rotations 
required to service  
all $n$ requests. 

Data locations on a disk may be 
specified by a radial coordinate $r$ which measures distance from the inner 
radius of the disk and by an angular coordinate $\theta $ which measures 
the angle with respect to some fixed ray. 
The seek function $f(\theta )$
is defined to be the radial distance which the disk 
head can travel, starting and 
finishing at rest in the radial direction, 
while the disk rotates through an angle 
$\theta $. The 
mechanics of disk head motion dictate that the function $f$ is convex. 
The choice of $f$ characterize a model of the physical disk. 
For a pair of disk locations $P,Q$ we let $d_f(P,Q)$ 
be the time required to move the disk head from location $P$ to 
location $Q$, with no radial velocity at the beginning and end of the movement.
$d_f$ always satisfies the triangle inequality but fails to be symmetric.
The disk scheduling problem is the traveling salesman problem on $n$ disk locations with 
distance function $d_f$. In general this problem is NP-Complete, \cite{ABZ}.
In this paper we assume that the seek function is linear, $f(\theta )=c\theta $. 
Since $c$ determines the 
function $f$ within the family of linear functions 
we will use $c$ to denote the function instead of $f$.  
While linearity of the seek function is a simplistic assumption many insights
can be gained from studying this case, some of which can be applied more generally. In the linear case Andrews, Bender and Zhang
(ABZ) have found a polynomial time algorithm for computing the optimal tour. Our goal is to provide average case analysis for
the length of the optimal tour and a better understanding of the ABZ algorithm. 

To that end we
consider a density function $p(\theta ,r )d\theta dr$ defined on the locations of data
on the disk which expresses the popularity of the data residing at location $(r,\theta )$. 
We let
$L_{n,c,p}$ denote the random variable which assigns to $n$ data locations $P_1,...,P_n$, sampled 
in accordance with the density function $p$, 
the length of the optimal tour between them w.r.t the distance function $d_c$.  
Given $c$ and $p$ the random variable $L_{n,c,p}$ concentrates w.h.p around it's mean value which asymptotically has the form 
$m_{c,p}\sqrt{n}$. We provide a geometrical description of $m_{c,p}$ in terms of space-time (Lorentzian) geometry.
A space-time (A compact 
Lorentzian manifold with boundary)
provides both a natural way of sampling points via the volume form and (at least locally) a natural partial order
which in relativity theory is interpreted as the past-future causal relation, \cite{BLMS}. We thus obtain a natural random poset structure, 
the causal structure
on $n$ points sampled w.r.t the volume form. such objects are known as causal sets in the (quantum) gravity literature, 
see \cite{So} for a survey. 
We show that the ABZ algorithm is intimately connected to a simple peeling process applied to causal sets
arising from the Lorentzian metric $ds^2=\frac{2}{c}p(\theta ,r)(dr^2-c^2d\theta ^2)$ on a disk drive.
Another algorithm (process) which can be interpreted in a similar fashion is patience sorting which computes the longest increasing subsequence (l.i.s) in a permutation. Patience sorting is a peeling process applied to causal sets associated with metrics of the form 
$ds^2=4p(x,y)dxdy$, where $p$ is a probability distribution on the unit square in the plane. The fact that increasing subsequences 
can be viewed as 
chains in a causal set was noted in \cite{BB}. As a result of the above observations we show that the optimal tour length and the length of a l.i.s can be estimated by the diameters of their associated space-time models. For the l.i.s problem this is essentially a reformulation of a result of Deuschel and Zeitouni, \cite{DZ}.
When $p$ depends only on $r$ we obtain the explicit formula $m_{c,p}=\sqrt{\frac{2}{c}}\int_0^1\sqrt{p(r)}dr$. To simplify matters we
first develop these computations without reference to space-time geometry.
 The Lorentzian geometry interpretation also leads to
further information on the behavior of the ABZ and patience sorting algorithms. In particular we generalize a result of Aldous and Diaconis 
on the statistics of pile sizes in patience sorting. The framework provided by the paper can also be used to analyze other discrete
processes such as airplane boarding, \cite{BBSS}, and the polynuclear growth model (PNG) in statistical physics, \cite{PS}.

We note that the analysis of disk scheduling (along with airplane boarding) 
seems to be the first application of Lorentzian geometry outside the domains of physics.

In the last section we take $p$ to be uniform. We can then use some refined estimates on the size 
of the longest increasing subsequence in random permutations to
obtain better bounds on the tour length.  

For the convenience of the reader we have also added an appendix containing a very brief introduction to Lorentzian geometry
and it's relativistic interpretation. 

Some of the results of the present
paper are based on the first part of \cite{Ba}. The Lorentzian interpretation and associated results are new to the present paper. \cite{Ba}
also contains some asymptotic results for more general convex seek functions. 

\section{The disk scheduling problem}
In this section we describe the batched disk scheduling problem
and briefly survey some results on increasing subsequences which are used
later on.
The batched disk scheduling problem  
was formally introduced in \cite{ABZ} section 2 which we follow
closely in the next subsection. 
\subsection{The batched disk scheduling problem}
A computer disk has the shape of an annulus, which geometrically can be thought of as a cylinder $C$.
For convenience we normalize the 
radial distance between the inner and outer circles to be 1. Each point on the
disk (cylinder) is represented by coordinates, $R=(\theta ,r )$ where
$0\leq r\leq 1$ is the radial distance from the inner 
circle and $0\leq \theta \leq 1 $ is the angle relative 
to an arbitrary but fixed ray. A complete circular angle is chosen to be 
of 1 unit instead of $2\pi $ for the convenience of future 
computations. The points $(0, r)$ and $(1,r)$ are thus identified. 
\footnote{ A disk drive actually consists of several annuli (platters). However this fact will have little significance to what follows} 

The disk drive has a recording head which is used for reading/writing data from/to the disk. The head is placed on an arm that can move radially in and out. 
In addition the disk rotates at a constant speed. The seek function $f(\theta )$ represents the maximum radial
distance that the disk head can travel starting and ending with no radial
motion, while the disk rotates through an angle $\theta $. 
Note that the angle $\theta $ is not limited to be less than 1 since moving a large 
radial distance may take more than one rotation of the disk. 
We will assume that $f$ is linear, $f(\theta )=c\theta $. 

Consider the infinite strip $U$ given by $0\leq r\leq 1$, 
$-\infty <t< \infty $.
We consider the operator $T$ on $U$ which is defined by 
$T(t,r )=(t+1,r)$. The cylinder $C$, representing the disk, is the quotient of $U$ 
with respect to the action of $T$. The quotient map $\pi : U\longrightarrow C$
is given by $\pi (t,r )=(t\; mod\; 1,r)$. We also let $I^2\subset U$ denote the set of points in $U$
with $0\leq t\leq 1$ which is the unit square. $I^2$ is a fundamental domain for the action of $T$ on $U$. 
We think of the $t$ coordinate as representing time, in units of complete disk revolutions. 
We normalize the time coordinate in such a way that the head of the disk is at angle 
$0$ at time $0$. It follows that $t \; Mod\; 1$ is the angle of the disk head at time $t$. 
If the disk head at time $t$ is in radial position $r$ 
then the disk position of the head is $\pi (t,r)\in C$.

Let $f$ be a convex function. 
We define a partial order $\leq_{f,hor}$ on $U$.
We say that $(t_1,r_1) \leq_{hor} (t_2,r_2)$ iff 
$t_2 \geq t_1$ and $f(t_2-t_1)\geq \vert r_2-r_1\vert $.
It is easy to verify that 
$\leq_{f,hor} $ is indeed a partial ordering.
If $f$ is the seek function of a disk we may interpret this relation as 
follows, $(t_1,r_1)\leq_{hor} (t_2,r_2)$ iff
a disk head which at time $t_1$ is in disk location 
$\pi (t_1,r_1)$ can reach at time $t_2$ the disk location $\pi (t_2,r_2)$.

Given a set of I/O requests to data locations, $\left\{R_1,...,R_n\right\}\subset C$, we formally add a request $R_0=(0,0)$ and denote the resulting subset of $C$ by $\bar{R}$. We also consider the set 
$\bar{R}^U=\pi^{-1}(\bar{R})$. 
 
A {\it tour} of $\bar{R}$ is a chain in the poset 
$(\bar{R}^U,\leq_{hor})$ starting at $(0,0)$
and ending at $(0,k)$, for some $k>0$, and containing one element from each set $\pi^{-1}(R_i)$
for all $i=1,...,n$. A tour can be thought of as giving the order and times in which the requests $R_i$ are to be serviced. 
The chain condition guarantees that the head of the disk can reach a request, from the location of the previous request, 
on time to provide service. We define the {\it service time} of the set of requests $\bar{R}$, denoted
$ST(\bar{R})$ to be equal to the minimal integer $k$ for which there 
a tour of 
$\bar{R}$ starting at $(0,0)$
and ending at $(0,k)$. Stated otherwise, the service time is the minimal number of disk rotations needed to service the 
requests in $\bar{R}$. 

Let $g$ be a concave function. 
We define a partial order $\leq_{g,ver}$ on $U$.
We say that  $(t_1,r_1)\leq_{g,ver} (t_2,r_2)$ iff  
$g(\vert t_2-t_1\vert )\leq \vert r_2-r_1\vert $. 

From now on we restrict our discussion to linear functions $f$ and omit $f$ (or $c$)
from subscripts. When considering a finite set $S$ in $U,C$ or $I^2$ we shall always assume that for any pair of points $(t_1,r_1), (t_2,r_2)$ we have $\vert f(\vert t_1-t_2\vert )\neq \vert r_1-r_2\vert $. 
By definition a set 
$S\subset U$ is a chain with respect to $\leq_{hor}$ iff it is an 
independent set with respect to $\leq_{ver}$ 
and likewise $S$ is an independent set with respect to $\leq_{hor}$ 
iff it is a chain with respect to $\leq_{ver}$. We shall refer to this property by saying that $\leq_{ver}$ and $\leq_{hor}$
are {\it complementary}. 

We define $\leq_{ver}$ and $\leq_{ver}$ on the unit square $I^2$ by restriction from $U$. 
We can also lift the relation $\leq_{ver}$ to $C$ by letting $R_1\leq_{ver} R_2$, $R_1,R_2\in C$ 
iff there exist $Q_1,Q_2\in U$
such that $\pi (Q_i)=R_i$ and $Q_1\leq_{ver} Q_2$. We note that for $Q_1,Q_2\in U$, such that $Q_1\leq_{ver}Q_2$ we have
$r(Q_1)< r(Q_2)$. It is easy to verify using this fact that $\leq_{ver}$ defines a partial order on $C$.
If $Q_1\leq_{ver} Q_2\leq_{ver} ...\leq_{ver}Q_n$ is a chain of $\leq_{ver}$ in $U$ then by definition the ordered set 
$R_i=\pi (Q_i)$ is a chain for $\leq_{ver}$ in $C$. Conversely, given a chain $R_1\leq_{ver}R_2\leq_{ver}...\leq_{ver}R_n$
in $C$ we can form a chain $Q_1\leq_{ver}...\leq_{ver} Q_n$ in $U$ such that $\pi (Q_i)=R_i$. We can construct 
the lifted chain as follows. Assume that $R_1,...,R_i$ have been lifted to $Q_1,...,Q_i$. By the definition of $\leq_{ver}$ in $C$ there are 
$Q'_i,Q'_{i+1}\in U$ such that $Q'_i\leq_{ver} Q'_{i+1}$ and such that $\pi (Q'_{i})=R_i$ and $\pi (Q'_{i+1})=R_{i+1}$.
Since $\pi (Q_i)=\pi (Q_i')$ we have an integer $k$ such that $Q_i=T^k(Q_i')$. Let $Q_{i+1}=T^k(Q_{i+1})'$ 
then by $T$ invariance
of the relation $\leq_{ver}$ in $U$ we obtain the desired extension of the lifted chain.

Given a finite partially ordered set $(S,>)$, there is a simple and well known
procedure which simultaneously finds a minimal decomposition of $S$ into independent sets and a maximal chain in $S$. 
Consider the set $S_1$
of minimal elements in $S$ with respect to $> $, let $S_i$ be the set of
minimal elements in $S-\cup_{j=1}^{i-1}S_j$. For each element of $s\in S_i$ 
construct pointers to all elements $s'\in S_{i-1}$ for which $s>s'$. Let
$S_k$ be the last nonempty set thus defined. We note that
the sets $S_i$ must be independent sets.
All maximal chains in $S$ are
obtained by following pointers from elements of $S_k$ all the way back 
to $S_1$. The length of the maximal chain is obviously $k$ which must therefore be the size of a minimal 
decomposition into independent sets. We note that the process can be applied to infinite posets as long
as the maximal length of a chain which ends at a given element $s\in S$ is bounded. 
The origin of this simple process may be traced back to G.Cantor's work on ordinal arithmetic.

\subsection{The optimal algorithm of Andrews, Bender and Zhang}
We consider an algorithm of
M.Andrews, M.Bender and L.Zhang for finding the optimal tour (up to a bounded additive constant)
for servicing a set of requests $\bar{R}$, \cite{ABZ}. 
We rephrase their algorithm in terms of the poset peeling process which was introduced earlier.

Let $\bar{R}\subset C$ be the set of requested data locations. 
Let $W_1,...,W_m$
be the sets obtained in the peeling process applied to the poset $(\bar{R}^U, \leq_{ver}$). Since
the $W_i$ are independent sets with respect to $\leq_{vert}$ they form chains
w.r.t $\leq_{hor}$. the chain structure provides a linear order on $W_i$. The number of independent sets
is by the properties of the peeling process  equal to the size of
of the maximal chain in $\bar{R}^U$ w.r.t $\leq_{ver}$. We denote the size of the maximal chain by 
$M_{U,ver}(\bar{R}^U)$ . 
By the chain lifting construction of the preceding subsection this equals
$M_{C,ver}(\bar{R})$, the maximal chain in $\bar{R}$ as a subset of $C$ w.r.t $\leq_{ver}$.   

Since the map $T$ is a partial order preserving isomorphism of $U$ onto itself,
the sets $W_i$ are invariant under $T$ and hence have the form 
$\pi^{-1}(S_i)$ for some $S_i\subset \bar{R}$. 

Given a chain of elements $v_1,..,v_l$ in a poset we will 
denote by $pred(v_i)$ the predecessor of $v_i$ in the chain, namely, 
the element $v_{i-1}$.

\

{\bf The ABZ Algorithm}: (see \cite{ABZ}) Start at the point $(0,0)\in U$.
proceed to the first point $w_{1,1}\in W_1$ for which $w_{1,1}\geq_{hor} (0,0)$.
Let $w_{1,1}, w_{1,2},...,w_{1,j(1)}=pred(T(w_{1,1}))$, be the elements of $\bar{R}^U$ in the chain $W_1$
between $w_{1,1}$ and $T(w_{1,1})$. Similarly, for general $1\leq i\leq m$, let $w_{i,1}$
be the minimal element $w\in W_i$ for which $w\geq_{hor} w_{i-1,j(i-1)}$.
Let $w_{i,1},...,w_{i,j(i)}=pred(T(w_{i,1}))$ be the elements of $\bar{R}^U$ in the chain $W_i$
between $w_{i,1}$ and $T(w_{i,1})$. Let
$k$ be the least integer for which $(0,k)\geq_{hor} w_{m,j(m)})$.
The output of the algorithm is the concatenated tour \newline
$\pi ((0,0)),\pi (w_{1,1}),...,\pi (w_{1,j(1)}), \pi (w_{2,1}),...,\pi (w_{m,j(m)}), \pi (0,k)$. 

\section{Estimating the optimal tour length}
\label{lis}
\subsection{A combinatorial estimate}
In order to analyze tour length for this algorithm 
it will be convenient to consider a slightly modified version of the ABZ 
algorithm.

Consider first the 
piecewise linear curves $L_i'$ in $U$ which join successive points of $W_i$ by
straight line segments. Since the $W_i$ are chains and the 
seek function $f$ is linear the curves
$L_i'$ form continuous chains with respect to $\leq_{hor}$. We also 
note that each element of $W_{i+1}$ lies above the line $L_i'$ 
since by definition an 
elements of $w_{i+1}\in W_{i+1}$ dominates an element of $w_i\in W_i$ 
w.r.t $\leq_{ver}$.
If the line $L_i'$ would have passed above we would get a 
point $l_i'$ of $L_i'$ with the same $t$ coordinate as $w_{i+1}$ which dominates $w_{i+1}$ 
and hence dominates $w_i$ w.r.t $\leq_{ver}$. 
By complementarity we conclude that there is a point $l_i'\in L_i'$ which is independent of $w_i\in W_i$ 
w.r.t $\leq_{hor}$ contradicting the chain condition. The $L_i'$ are also invariant under $T$ since $W_i$ is $T$ invariant.
Finally we note that the piecewise linear curves $L_i'$
can be considered as graphs of functions on $t$. Let $L_1=L_1'$ and 
define the curves $L_i$ inductively 
to be the graph of the maximum of the functions $L_i'(t)$ and $L_{i-1}(t)$.
Since by the argument above points of $W_i$ lie above $L_{i-1}$ and by definition lie on 
$L_i'$ they also lie on $L_i$. By definition the curves $L_i$  
form chains w.r.t $\leq_{hor}$ and at  
Every point $L_{i+1}$ lies above or coincides with $L_i$. $L_i$ are again $T$ invariant by definition.

Finally consider the line segment $J$ defined by $r=ct$, passing from the point $(0,0)$ to the point $(1/c,1)$. 
By definition $J$ forms a chain w.r.t $\leq_{hor}$. Let $l_{i,1}$ be the first meeting point of $J$ 
with $L_i$. Since $J$ begins at the point $(0,0)$ and $L_{i+1}$ is above $L_i$
we see that $l_{i+1,1}$ is either further along $J$ or equal to $l_{i,1}$. We let $(0,0)$ be $l_0$
  
We can now state the modified algorithm

\

{\bf Modified ABZ algorithm}:

Let $l_{i,1},l_{i,2},...,l_{i,h(i)}=T(l_{i,1})$ be the points of $\bar{R}^U$ on $L_i$ between $l_{i,1}$ and $T(l_{i,1})$, 
inclusive. The modified ABZ tour consists of the concatenation \newline 
$(0,0), l_{1,1},...,T(l_{1,1}),T(l_{2,1}),...,T^{i-1}(l_{i,1}),...,T^i(l_{i,1}),...,T^m(l_{m,1}), (k,0)$, 
where $k$ is the smallest integer such that $T^m(l_{m,1})\leq_{hor} (k,0)$. 

\begin{theo} 
\label{combinatorial}
$M_{C,ver}(\bar{R})-1-1/c \leq ST(\bar{R})\leq M_{C,ver}(\bar{R})+1+2/c$
\end{theo}
{\bf Proof }: We consider the upper bound first. We have $t(T^m(l_{m,1}))=m+t(l_{m,1})\leq m+1/c$. 
The point $q=(t(T^m(l_{m,1}))+cr(T^m(l_{m,1}),0)$ satisfies $q\geq_{hor}T^m(l_{m,1})$, hence \newline 
$k=[t(T^m(l_{m,1}))+(1/c)r(T^m(l_{m,1})]+1\leq m+1/c+1/c+1=m+2/c+1$.

We now consider the lower bound. Let $Q_1,Q_2\in U$ such that $Q_2\geq_{ver} Q_1$. Let $Q_2'\in U$ be such that 
$\pi (Q_2')=\pi (Q_2)$ and $Q_2'\geq_{hor} Q_1$. If $t(Q_2)\geq t(Q_1)$ then
by definition of the relation $\leq_{ver}$ we have $(t,r(Q_2))\geq_{ver} Q_1$ for all $t(Q_1)\leq t \leq t(Q_2)$.
We conclude by complementarity that $t(Q_2')\geq t(Q_2)+1=t(Q_1)+(t(Q_2)-t(Q_1))+1$. If $t(Q_2)\leq t(Q_1)$ or 
$Q_1\geq_{ver} Q_2$ then
obviously the same inequality holds. 
Let $R_{i_1},...,R_{i_m}$ be a maximal chain w.r.t $\leq_{ver}$ among the
requests $\bar{R}$. Let $Q_1,...,Q_m\in U$ be a lifting w.r.t $\pi $ of the maximal chain to a maximal chain w.r.t $\leq_{ver}$
on $U$. Consider an optimal tour in $U$. Assume that the points $Q_i$ are indexed by the order in which the optimal tour visits
their equivalence classes w.r.t $T$, this does not coincide necessarily with their order in the chain. Consider $Q_1$ and $Q_2$
let $Q_i'$ be the point in the equivalence class of $Q_i$ which the tour visits. 
By the chain condition we have either $Q_1\geq_{ver}Q_2$ or vice versa, in any case we conclude that 
$t(Q_2')\geq t(Q_1)+(t(Q_2)-t(Q_1))+1$. More generally, assume that $Q_i'=T^{l_i}(Q_i)$. Applying the argument above to
$Q_i'$ and $T^{l_i}(Q_{i+1})$ we conclude that $t(Q_{i+1}')\geq t(Q_i')+t(T^{l_i}(Q_{i+1}))-t(Q_i'))+1=t(Q_i')+(t(Q_{i+1})-t(Q_i))+1$,
which leads to $t(Q_m')-t(Q_1')\geq (m-1)+\sum_{i=1}^{m-1}t(Q_{i+1})-t(Q_i)=(m-1)+(t(Q_m)-t(Q_1))$. Since $Q_1$ and $Q_m$
are comparable w.r.t $\leq_{ver}$ we have $t(Q_m)-t(Q_1)\geq -1/c$ which leads to the desired estimate. {\em q.e.d}  
\subsection{The probabilistic setting}
We wish to provide average case analysis for the service time
$ST(\bar{R})$ as the number of requests, $n$, goes to infinity.
We assume that we are given a probability distribution 
$\mu_p=p(\theta ,r )d\theta dr$, where we will assume that $p$ is smooth,
on disk locations which reflects the popularity of the data stored in various disk locations. 
We assume that the $n$ requested locations are sampled 
independently from the distribution $\mu_p$. We may also consider a dependent sequence of requests in which case
$p$ will denote the asymptotic distribution of requests, assuming one exists.
We will be interested in asymptotic behavior of $ST(\bar{R})$ as $n$ 
tends to infinity. 
In this context the phrase ``with high probability'', w.h.p for short,
refers to an event that occurs with probability approaching 1 
on the probability space of all sets of $n$ requests as $n$ approaches
infinity.

{\bf Remark}: We remark that in practice $p$ will depend only on $r$ and not on 
$\theta $. This is due to the way data is laid out on a disk. 
Logically consecutive blocks of data are laid out in tracks for which $r$ is fixed. 
To produce bias in the angular direction $\theta $ an application would have to request 
blocks which are in arithmetic progressions with jump $d$ equal to a multiple of the number of blocks in a disk track. 
Applications are usually unaware of the physical layout of data on disks, It is therefore very improbable that an application would randomly or even intentionally be 
able to produce bias in the $\theta $ direction. 
On the other hand radius dependent bias occurs very frequently. Empty regions of the disk 
are of the form $r_1\leq r\leq r_2$ because of data layout, also,
files occupy regions of this form resulting in bias due to differing file popularities. 
Nonetheless, it is instructive to work out the case of a general function $p(\theta ,r)$.  
\subsection{Increasing subsequences} 
We say that a sequence of points $z_1=(x_1,y_1),...,z_k=(x_k,y_k)$ in the 
plane is {\it increasing} if $x_i\geq x_j$ iff $y_i\geq y_j$ for all 
$1\leq i,j\leq k$. 

We can define a partial order $\geq_{inc}$ on points of the plane by saying
that $(x_1,y_1)\geq_{inc} (x_2,y_2)$ iff $x_1\geq x_2$ and $y_1\geq y_2$.
An increasing subsequence is by definition a chain in this partially ordered set.

Let $\mu_q=q(x,y)dxdy$ be a distribution on a square with vertices $(0,0), (0,a), (a,0), (a,a)$ with bounded density function $q$. The following
theorem of Deuschel and Zeitouni computes the asymptotics of
the longest increasing subsequence among $n$ points in the unit square sampled
w.r.t $\mu_q$.
   
\begin{theo} \label{DZ} (\cite{DZ})

Let $S$ be a set of $n$ points in the unit square 
chosen with respect to $\mu_q$. 
Denote by $K$ the size of the largest increasing subset of $S$ w.r.t $\leq_{inc}$. Then  

\begin{itemize} 
\item[1) ] For all $\varepsilon >0$, w.h.p, $\vert K-\ell_{max} \sqrt{n}\vert <\varepsilon \sqrt{n}$. Here $\ell_{max}$ 
is given by 
$\ell_{max}=Max_{\phi } \ell (\phi )$, where $\ell (\phi )$ is the functional 
\begin{equation}
\label{length}
\ell (\phi )=2\int_0^a\sqrt{\phi '(x)q(x,\phi (x))}dx 
\end{equation} 
and $\phi $ runs through all differentiable 
nondecreasing functions on the unit
interval with boundary conditions $\phi (0)=0$ and $\phi (a)=a$.

\item[2) ] For any $\varepsilon ,\delta >0$, 
w.h.p, an increasing subset of size $(K-\varepsilon )\sqrt{n}$ can
be found in a $\delta $ neighborhood of $\phi $ if $\phi $ maximizes
the functional $\ell (\phi )$. Here a $\delta $ neighborhood refers to all points
which are at a distance less then $\delta $ from a point of the form
$(x,\phi (x))$.
\end{itemize}
\end{theo}

\subsection{Estimating $M_{C,ver}(\bar{R})$}
Assume that a set of $n$ I/O requests, $\bar{R}=\left\{R_1,...,R_n\right\}$ are chosen w.r.t the location 
distribution given by $p(\theta, r)d\theta dr$. 
Our goal is to study the asymptotics of the random variable 
$ST(\bar{R})$. By theorem \ref{combinatorial} we may equivalently study the asymptotic behavior of
the random variable $M_{C,ver}(\bar{R})$. We do so using 
theorem \ref{DZ}. We define the pull back $p_U$ of the function $p(\theta, r)$ to $U$ by letting  
$p_U(t,r))=p(\pi (t,r))$ for $(t,r)\in U$. We recall that the seek function $f$ is given by $f(t)=ct$.
\begin{theo}
\label{estimate}
1) w.h.p $\vert M_{C,ver}(\bar{R})-l(c,p)\sqrt{n}\vert< \varepsilon \sqrt{n}$, where $l(c,p)$ is given by
$Max_{\psi }l (\psi )$ and 
\begin{equation}
\label{lengthC}
l(\psi )=\sqrt{\frac{2}{c}}\int_0^1\sqrt{p_U(\psi (r), r)(1-c^2\psi '(r)^2)}dr
\end{equation}
and the maximum is taken over all differentiable functions $\psi (r)$ with $\vert \psi '(r)\vert \leq 1/c$.
  
\

2) If $p=p(r)$ then $l(c,p)=\sqrt{\frac{2}{c}}\int_0^1\sqrt{p(r)}dr$
\end{theo}
{\bf Proof}: We prove the theorem first for the case of $c=1$. Define the rectangles $D_{i,j}$ $i=0,1,2,3$, $j=-4,...,4$
in $U$ consisting of the points
with coordinates $i/4\leq r\leq (i+1)/4$ and $t\; Mod \; 1$ with $j/4\leq t \leq (j/4)+1 $. Given a chain $M$ in $C$ w.r.t $\leq_{ver}$ 
we can find a function $g=g_M$ $g:\left\{ 0,1,2,3\right\}\rightarrow \left\{ -4,...4\right\}$ and a lifting of $M$ as constructed in the previous section 
which is a chain in 
$D_g=\cup_i\;  D_{i,g_M(i)}$, $i=0,1,2,3$. This is the case since the $\theta $ coordinate of
elements in $M$ can vary by at most $1/4$ as $r$ varies by $1/4$ and any $\theta $ interval of size $1/4$ thereby missing at least one vertical 
line of the form $j/4$.
We conclude that the length of the maximal chain is the maximum,
over all functions $g$ as above, of the maximal size chain in $D_g$. We further note that each $D_g$ is a fundamental domain for the translation 
action $T$ on $U$ and hence, the probability distributions $p(\theta, r)d\theta dr$ lifts to a distribution on $D_g$. 
Fix a function $g$ and a (piecewise) 
differentiable curve
$(r,\phi (r))$ in $D_g$ with $\vert \phi '\vert \leq 1$. We consider the size of a maximal chain in $D_g$ in an $\varepsilon $ neighborhood 
of the graph of $\phi $. It is easy to see that any chain will be in the $\varepsilon$ neighborhood of such a graph.

Let $S$ be the linear transformation which rotates the plane clockwise by 45 degrees. It
is easy to verify that for points $R_1,R_2$ in $D_g$ we have $R_1\leq_{ver} R_2$ iff $S(R_1)\leq_{inc} S(R_2)$. 
We can now apply theorem \ref{DZ} to $S(D_g)$ with the induced distribution. Up to a translation, $S(D_g)$
will be contained in a square as in theorem \ref{DZ}. 
We observe that curves $(x, \phi (x)) $ with $\phi ' (x)\geq 0$ are mapped by $S^{-1}$ to graphs of functions $\psi (r)$
which satisfy $\vert \psi '(r) \vert \leq 1$. At the level of differentials we have along the graph of $\phi $, the formula 
$dy=\phi '(x)dx$. Also $dx=(dr+d\theta )/\sqrt{2}$ and similarly 
$dy=(dr-d\theta )/\sqrt{2}$. We conclude that  
$dxdy=\phi ' (x)(dx)^2=((dr)^2-(d\theta )^2)/2=((dr)^2(1-(\psi '(r))^2)/2$. We see that the functional whose maximal value is $l(c,p)$ coincides with the functional $\ell (\phi )$ of theorem \ref{DZ} as written in the $\theta, r$ coordinate frame 
via the mapping $S^{-1}$. Since $\leq_{ver}$ 
and $\leq_{inc}$ are identified via $T$ we obtain the theorem for $\leq_{ver}$ restricted to $D_g$ and since any chain
of $\leq_{ver}$ in $C$ can be lifted to one of the $D_g$ and the number of functions $g$ is finite we obtain the theorem.  

For general values of $c$, the transformation $V_c(\theta, r)=(c\theta, r)$ transforms $\leq_{c,ver}$ to $\leq_{1,ver}$ and 
we can apply the preceding arguments. The factor $\frac{1}{c}$ comes from the effect of $V_c$ on the probability distribution.
  
Let $p=p(r)$ be a density which only depends on $r$. Let $0=r_0<r_1<r_2<...<r_{n-1}<r_n=1$ and let $d_i$ 
$i=1,...,n$ be arbitrary non negative coefficients with 
$\sum_{i=1}^n d_i(r_i-r_{i-1})=1$. Let $q=q(r)$ be a density function of the form  
$\sum_{i=1}^n c_i\chi_{I_{r_{i-1},r_i}}$ where $\chi_{I_{r_{i-1},r_i}}$ is the
characteristic function of the interval $[r_{i-1},r_i]$. We will assume that the set $r_0,...,r_n$ contains
the points $0,1/4,1/2,3/4,1$.

Since it is known that finite linear combinations of characteristic 
functions of intervals are dense in the space of Lebesgue measurable functions
it is enough to prove the theorem for such densities.

Consider the set $J_i$  
consisting of points $(t,r)$ with $r_{i-1}\leq r\leq r_i$. Fix a function $g$ as above and consider $J_{i,g}=J_i\cap D_g$. Since we have 
included the points $k/4$, $k=0,...,4$ in the set $r_j$ we have for any 
$g_1,g_2,i$ a translation in $U$, $T_{g_1,g_2,i}:(t ,r)\longrightarrow 
(t+t(g_1,g_2,i),r)$ which maps $J_{i,g_1}$ to $J_{i,g_2}$. We note that translations preserve the functional in our theorem. The image of
$J_{i,g}$ under $S$ is a rectangle rotated by $\pi /4$. We can assume that the bottom vertex is the origin by
shifting the image. The vertices of the rotated rectangle are then of the form $(0,0), (-a,a), (b,b), (b-a,b+a)$ for
some positive $a,b$. It is shown in \cite{DZ} that given a constant density function on a convex set $A$ and
$u=(x_1,y_1),w=(x_2,y_2)\in A$, the straight line between $u$ and $w$ is the unique curve which maximizes the
functional \ref{length}, subject to the initial conditions $\phi (x_1)=y_1$ and $\phi (x_2)=y_2$.
Since $q$ is constant
on $J_{i,g}$ we see that the curve $\phi $ maximizing the functional on $T(J_{i,g})$ is a straight line between two boundary points
$u$ and $w$. The value of the functional $\ell_{max}=\ell (\phi )$ is then given by $2\sqrt{q(x_2-x_1)(y_2-y_1)}$. We also notice that
$u$ must lie on the interval $A=[(0,0),(-a,a)]$ and $w$ on the parallel interval $B=[(b,b),(b-a,b+a)]$. 
Since both intervals have a slope of $-1$ we notice that in fact $(x_2-x_1)+(y_2-y_1)$ is independent of $u\in A$ and 
$w\in B$. The product is then maximized when $x_2-x_1=y_2-y_1$ in which case the line $\phi $ satisfies $\phi '=1$.
Pulling back the result we see that any of the curves $\psi (r)=constant$ in $J_{i,g}$ maximizes the functional among all
curves in $J_{i,g}$ and these are the only maximizing curves. Consider a vertical line of the form 
$\psi (r)=constant$ and
an arbitrary non vertical curve $\hat{\psi } $. 
Let $g_{\psi }$ and $g_{\hat{\psi } }$ be such that $\psi $ is contained in $D_{g(\psi )}$ and $\hat{\psi } $ is contained in
$D_{g(\hat{\psi })}$. Let $\psi_i=\psi \cap J_{i,g(\psi )}$ and $\hat{\psi }_i=\hat{ \psi } \cap J_{i,g(\hat{\psi })}$.
The line $T_{g(\psi ),g(\hat{\psi } ,i}(\psi )$ is a vertical line in $J_{i,g(\hat{\psi })}$ and hence 
$l(\psi_i)=l(T_{g(\psi ),g(\hat{ \psi }),i}(\psi_i ))\geq l(\hat{\psi } )$. Since $\hat{\psi } $ is not vertical there exists an $i_0$
for which the inequality is strict. Since $l(\psi )=\sum_i l(\psi_i)$ and $l(\hat{\psi } )=\sum_i l(\hat{\psi }_i)$ we see
that $l(\psi )>l(\hat{\psi })$. Since this inequality holds for any density step function density $q$ it holds for
$p(r)$. Computing $l$ on an arbitrary vertical line yields the explicit formula.
{\em q.e.d}

\

As a consequence of the proof of theorem \ref{estimate} we obtain the following corollary.
\begin{coro}
Let $\bar{R}\subset C$ be a set of $N$ requests. The modified ABZ tour on $\bar{R}$ can be computed 
in time $O(n\log (n))$. It provides a tour which is optimal within an additive constant.
\end{coro}

\noindent {\bf Proof}: It is sufficient to establish the running time. By rescaling we can assume that $c=1$.
 The main computational step in the modified ABZ algorithm is the computation 
of the peeling process decomposition $W_1,...,W_m$ for the set $\bar{R}^U$ w.r.t $\leq_{ver}$. 
When $c=1$, points in $I^2$ can dominate or be dominated only by points in the set $B\subset U$ satisfying $-1\leq t\leq 1$. Consequently, 
 the restriction to $I^2$ of 
the peeling process applied to points in $B$ coincides with the restriction to $I^2$ of the peeling process applied to all $U$. We have noted in the proof of theorem \ref{estimate} that by a 45 degree rotation the relation $\leq_{ver}$
translates into the relation $\leq_{inc}$. It is therefore sufficient to apply the peeling process to the set of size $3N$, $S(\bar{R}^U\cap B)$ w.r.t $\leq_{inc}$. It is well known that patience sorting, \cite{AD}, provides an 
$O(nlog(n))$ algorithm for peeling a subset in the plane w.r.t $\leq_{inc}$. {\em q.e.d.}
 
\section{Lorentzian geometry and disk scheduling}
\label{spacetime}
In this section we reformulate and extend the results of the previous section in terms of Lorentzian geometry.
For the convenience of the reader we have added a very brief introduction to Lorentzian geometry in the appendix.
We will use some basic results on Lorentzian geometry in this section. These results can be found in \cite{Pe}, especially
in chapter 7. 
We consider Lorentzian metrics on $C$ the cylinder, $U$ the infinite strip and $I^2$ the unit square. 
We recall the mapping 
$\pi:U\longrightarrow C$, $\pi(t,r)=(t\; Mod\; 1,r)$. We consider
the standard coordinates on them $(\theta ,r)$, $(t,r)$ and $(x,y)$.

\footnote{Strictly speaking $(\theta ,r)$ is not a coordinate 
system on $C$ due to the identification of $(0,r)$ with $(1,r)$. 
We can define coordinate patches on $C$ similar to the subsets $D_f$
in the proof of theorem \ref{estimate}. The transition functions between the coordinates 
on different patches are given by translations 
(rotations in $C$) which do not affect the form of the metric}

Let $p$ be a density function on either $C$ or $I^2$. 
We will define metrics by their associated forms $ds^2$. Given $p$ we define the metric $g_{C,ver}$ by
$ds^2=\frac{2}{c}p(\theta ,r)(dr^2-c^2d\theta ^2)$. We define the metric $g_{U,ver}$ on $U$ 
by $ds^2=\frac{2}{c}p(\pi (t,r))(dr^2-c^2dt^2)$. We define
$g_{U,hor}=-g_{U,ver}$. On $I^2$ we consider the metric $g_{I^2,inc}$ 
given by $ds^2=4p(x,y)dxdy$. For $g_{C,ver}$ and $g_{U,ver}$
we define the future
pointing vector to be $(0,1)$ globally. For $g_{U,hor}$ we define 
the future pointing vector to be $(1,0)$ at all points and for 
$g_{I^2,inc}$ we define the future pointing vector to be $(1,1)$ at all points. 
With these definitions it follows immediately
that the partial orders we considered coincide with the past-future causal relation induced by the corresponding 
Lorentzian metrics.

{\bf Remark}: We can also define on $C$ the metric $g_{C,hor}=-g_{C,ver}$. 
The metric $g_{C,hor}$ does not have an associated 
partial order (causal structure) since there are time-like closed curves, 
in fact any horizontal circle of the form $r=constant$ is a closed time-like curve for this metric.   

Given a point $z=(\theta ,r)$ we let $h(z)$ be the maximal length of a curve ending at $z$. A maximal length curve ending at $z$ 
exists by corollary 7.7 of \cite{Pe}. We will refer to $h(z)$ as the {\it height} of the point $z$.
Let $C_{\tau }$ denote the curve consisting of points $z$ such that $h(z)=\tau $ in the space-time $(C,g_{C,ver})$. 
The following result provides a reformulation and extension of theorem \ref{estimate} in terms of
Lorentzian geometry.
\begin{theo}
\label{disk}
Let $p(\theta, r)$ be a request probability distribution 
on a disk $C$. Let $g$ be the Lorentzian metric induced by
$ds^2=\frac{2}{c}p(\theta ,r)(d\theta ^2-c^2r^2)$. Let $diam(C)$ denote the maximal length of a time-like curve in this model. 
Let $\bar{R}$ be a set of $n$ requests chosen w.r.t to the density distribution $p$ and let
$ST(\bar{R})$ be the number of disk rotations needed to service all 
the requests in $\bar{R}$ with an optimal policy, then:

\

1) For all $\varepsilon >0$, w.h.p 
\begin{equation}
\label{A}
\vert\frac{ST(\bar{R})-diam(C)}{\sqrt(n)}\vert <\varepsilon \sqrt{n} 
\end{equation}

\

2) When $p=p(r)$ we have $diam(C)=\sqrt{\frac{2}{c}}\int_0^1\sqrt{p(r)}dr$. 

\

3) Let $n(\tau )$ be the number of I/O which were serviced in the first $\tau \sqrt{n}$ disk rotations.
For all $\varepsilon >0$ we have w.h.p 
\begin{equation}
\label{B}
\vert \frac{n(\tau )-\frac{1}{2}\int_0^{\tau }\ell (C_h)dh}{n}\vert <\varepsilon
\end{equation}
where $\ell (C_h)$ is the length of $C_h$ w.r.t $-g$.
\end{theo}

\noindent {\bf Proof}: Given the definitions, parts 1 and 2 are simply a reformulation
of theorem \ref{estimate}. We therefore need to prove part 3. Fix some angle $\theta $. The function $h_{\theta }(r)=h((\theta ,r))$
is a strictly increasing function. Given some $\tau $ we conclude that either there is a unique point of the form $\theta ,g_{\tau }(\theta )$
of height $\tau $, or $h((\theta ,r))<\tau $ for all $0\leq r\leq 1$. We define $g_{\tau }(\theta )$ to be 1 in the latter case
and define $\tilde{C}_{\tau }$ to be the graph of $g_{\tau }$. The curve $C_h$ must be space-like by definition and hence it is a causal curve for the metric $-g$, see \cite{Pe} page 17 for a definition. It is also continuous at interior points 
of $C$. As a causal curve w.r.t $-g$ it also has a finite length 
w.r.t $-g$ according to \cite{Pe} page 54. We also know that $C_h$ is differentiable almost everywhere, \cite{Pe} page 17.
Let $X(\tau )$ be the region bounded by $C_0$ 
and $C_{\tau }$. By theorems \ref{DZ}, \ref{estimate} and their proofs we know that for every $\delta >0$, w.h.p, the I/O
which were serviced during rotation $\tau \sqrt{n}$ lie in a $\delta $ neighborhood of $C_{\tau }$. Letting
$\delta $ tend to zero we see that for each $\varepsilon >0$ w.h.p 
$1-\varepsilon <\frac{n(\tau )-(\int_{X_{\tau}}p(\theta ,r)d\theta dr)n}{n}<1+\varepsilon $. By equation (\ref{volume}) the Lorentzian
area element is $2p(\theta ,r)d\theta dr$ hence $\int_{X_{\tau }}p(\theta ,r)d\theta dr=1/2 vol(X_{\tau })$. It remains to show 
that $vol(X_{\tau })=\int_0^{\tau }l(C_h)dh$. 
For simplicity we shall presently assume that $g_h(\theta )<1$ for all $\theta $, that is the curve $C_h$ does not intersect the upper boundary curve $r=1$.

At this point we note that the $C_h$ are the circles of radius $h$ (from $C_0$), 
$X_{\tau }$ is the ball of radius $\tau $ and hence the desired formula is the classical formula of computing the volume (area) 
of a ball
from the volumes (circumferences) of spheres (circles), the proof of which follows from the fact that circles are perpendicular to (the geodetic) radial curves. 

We make this argument more precise. Let $z_1,z_2\in C_{\tau }$, not necessarily different, and consider maximal curves 
$\gamma_1, \gamma_2$ ending at $z_1,z_2$ respectively.  
By proposition 7.8 of \cite{Pe} $\gamma_i$ are geodesics.
The $\gamma_i$ cannot meet at any point. assume to the contrary that they meet at a point $y$. The
portions of $\gamma_i$ from $C_0$ to $y$ must also be maximal and therefore of the same length. The $\gamma_i $, 
being distinct geodesics must also have distinct 
tangent directions at $y$. Following $\gamma_1$ from $C_0$ to $y$ and then $\gamma_2$ from $y$ to $z_2$ would yield a maximal 
curve ending in $z_2$ which is not differentiable at $y$ and hence not a geodesic, see also \cite{Pe} page 57  
for this type of argument. By maximality, the $\gamma_i$ (being radial curves) 
must also be perpendicular to all $C_h$, $h< \tau$ at points in which they are differentiable, see \cite{Pe} theorem 7.27(b). 
Consider a point 
$z\in C_h$ with more than one maximal curve. We shall call such a point a {\it bad } point. We denote by $Z$ 
the set of all bad points. A point which is not bad will be called {\it good}.
We may parameterize the past pointing geodesics ending at $z\in Z$ by their tangent angle $\phi $
at $z$. The set of angles $\phi $ corresponding to maximal curves is closed since by upper semi continuity of the length function, \cite{Pe} theorem 7.5, a limiting curve of maximal curves will also be maximal. We conclude that there are geodesic maximal curves
$\gamma_{1,z},\gamma_{2,z}$ with endpoint $z$ with maximal and minimal tangent angles $\phi_1<\phi_2$. Given a bad $z$ we may consider 
all points in $C$ which lie in the curved triangle whose boundary curves are $C_0$, $\gamma_{1,z}$ and $\gamma_{2,z}$. denote 
this triangle by $T_z$ and let $T=\cup_{z\in Z}T_z$. Since maximal curves do not intersect, the points of $T_z$, either do not lie on any maximal curve 
or, can only lie on maximal curves ending at $z$. On the other hand a point of $C_{\tau }$, $\tau <h$, outside of $T$ will lie on a unique maximal curve ending at a good point. To see this, assume that $u\in C_{\tau }$ is such a point and assume that no maximal curve passes through $u$. Let $u_1,u_2\in C_{\tau }$ be the closest points to $u$ in $C_{\tau }$ through which maximal curves pass. As argued before, 
such points exist since limiting curves of maximal curves are maximal. The maximal curves through $u_1,u_2$ will end at points $z_1,z_2$. $z_1\neq z_2$ since otherwise $u\in T_{z_1}$ by definition. Any maximal curve of any point in $C_h$ between $z_1$ and $z_2$
will meet $C_{\tau }$ in a point $u_3$ in the interior of the interval of $C_{\tau }$ between $u_1$ and $u_2$ which contains $u$,
resulting in a contradiction to the choice of $u_1,u_2$. For any bad $z$ and $\tau <h$
we consider the interval $I_{\tau ,z}=T_z\cap C_{\tau }$. Since maximal curves do not intersect in $C_{\tau }$, $\tau <h$
we conclude that for any pair of bad points $z_1\neq z_2$ and any $\tau <h$ $I_{\tau ,z_1}$ and $I_{\tau z_2}$ are disjoint.
In particular the set of bad $z$ is countable. Let $B_{\tau } =\sum_{z\in Z} \ell (I_{\tau ,z})$. Since the intervals $I_{\tau ,z}$
are disjoint we have $B_{\tau }\leq \ell (C_{\tau })$. Consider a point $u\in I_{\tau ,z}$. The distance between $u$ and $z$ is at most
$h-\tau $ since otherwise there would be a curve from $C_0$ to $z$ passing $u$ of length greater than $h$. Let $T^{dh}$ be the set
of point in $T$ of height $h-dh\leq \tau \leq h$. We wish to show that $vol(T^{dh})=o(dh)$. 
For a given bad $z$ we consider
$T_z^{dh}=T_z\cap T^{dh}$. We consider the Lorentzian metric at a small neighborhood of $z$. 
We assume that the neighborhood is chosen to be simply convex, \cite{Pe} pages 5-6, so there are unique geodesics between any pair of 
points in the neighborhood. 
such neighborhoods always exist, see \cite{Pe} page 5. We consider Riemann normal coordinates, see \cite{Pe} pages 5-6 and \cite{Ca} pages 112-113, on a small enough neighborhood of $z$.  The metric in Riemann normal coordinates is constant up to second order terms, thus, we can compute lengths 
and areas up to a multiplicative error of $1+dh$ assuming a constant metric. 

In these coordinates the geodesics emanating from $z$ are given by straight lines. Also circles around $z$ coincide with the circles in Minkowski space. 
The lines $\gamma_{1,z}$ and $\gamma_{2,z}$ meet the curve $C_{h-dh}$
at a distance of $dh$ from $z$. $T_z^{dh}$ is contained in the sector $U_z$, bounded by $\gamma_{1,z}, \gamma_{2,z}$ and the
circle $S_{z,dh}$ of radius $dh$ around $z$. Consider $v_i$, $i=1,2$, the meeting points of $\gamma_{i,z} $ and $S_{z,dh}$. In the constant metric approximation 
The length of the geodesic between the $v_i$ grows linearly in $dh$
and can be written as $c_zdh$, where $c_z$ depends on the angle $\phi_{2,z}-\phi_{1,z}$ and the density $p(z)$. 
We consider the circle section $S_z$ of radius $dh$ between the $v_i$. Since the
the geodesic is the maximal length curve between the meeting points we have $c_zdh\geq \ell (S_z)$. By a simple calculation 
$vol(U_z)=\int_{\rho =0}^dh\ell S_{z,\rho }d\rho \leq \ell (S_{z,dh})dh \leq c_zdh^2$ in the constant metric approximation. We conclude
that $vol(T_z)\leq c_zdh^2 (1+dh)$. Notice however that replacing the intervals $I_{h-dh ,z}$ by the geodesics connecting their endpoints
still leads to a causal curve, hence of finite length, $l$. Since $\sum_zc_zdh+O(dh^2)\leq l$ we conclude that $\sum_zc_z$ is finite and 
$vol (T^dh)\leq (\sum_zc_z)dh^2(1+dh)=
O(dh^2)=o(dh)$. Let $X^{dh}=X_h-X_{h-dh}$. On $X^{dh}-T^{dh}$ we can define a vector field by the tangent of the unique geodesic passing through a point.
This vector field is almost everywhere orthogonal to the level curves $C_{\tau }$ and in particular to $C_h$. the area of $X^{dh}-T^{dh}$ 
is then $\ell (C_h)dh(1+dh)$. Combining with the estimate for $T^{dh}$ proves the theorem. 
We now consider the case where the curve $C_h$ does meet the upper boundary curve $F$ given by $r=1$. We define $F_{h,dh}$
to consist of the points on $z\in F$ of height $h-dh\leq h(z)\leq h$. We can apply the preceding arguments with 
$C_h$ replaced by the causal curve $\tilde{C}_{h,dh}=C_h\cup (F_{h,dh})$. The curve $C_h$ is the intersection of the monotone family of curves $\tilde{C}_{h,dh}$. Since the length of a portion of $F$ is bounded up to a constant multiplicative 
factor by it's measure we see $\ell (\tilde{C}_{h,dh})$ converges to $\ell (C_h)$ as $dh$ tends to 0 completing the proof.
{\em q.e.d}

Using essentially the same proof
we obtain the following more general version which is stated entirely in Lorentzian geometric terms 
\begin{coro}
\label{lorentzformulation}
Let $C$ be a compact subset of a strongly causal Lorentzian manifold, which satisfies $vol(C)=2$. Let $<$ denote the causal partial order induced on $C$. Consider the density distribution $p=vol/2$ on $C$. For any pair of points $x,y\in C$ with $x<y$ let $d(x,y)$
be the length of the maximal curve in $C$ beginning at $x$ and ending at $y$. Let $L_{x,y,n}$ denote the random variable which
assigns to $n$ points in $C$ sampled according to $p$ the length of the longest chain of the form $x<z_1<...<z_k<y$.
 The following statements hold:

\

1) for any pair of points $x,y\in C$ with $x<y$, w.h.p $1-\varepsilon <\frac{L_{x,y,n}}{d(x,y)\sqrt{n}}<1+\varepsilon $ 

\

2) Let $n_{\tau }$ denote the number of points in the first $\tau \sqrt{n}$ layers of the peeling process applied to $n$
sampled points, then for all $\tau >0$, 
$1-\varepsilon < \frac{n_{\tau }}{\frac{1}{2}(\int_0^{\tau }\ell (C_h)dh)n}<1+\varepsilon $.
\end{coro}

\

{\bf Example}: Let us consider the second statement with $C$ being the unit square and $ds^2=4dxdy$. In this case the peeling process
coincides with patience sorting, the layers corresponding to card piles, see \cite{AD}. The curve $C_h$ consists of the points satisfying $2\sqrt{xy}=h$ or $y=\frac{h^2}{4x}$. We use $x$ to parametrize the curve and $dy=\frac{-h^2}{4x^2})$ hence the curve length is 
$\int_{\frac{h^2}{4}}^1 \frac{h}{x}dx=2hln (\frac{2}{h})$. 
Part 2 then states that asymptotically $n_{\tau }=\int_0^{\tau }hln(\frac{2}{h})dh$
thus reproving theorem 12 of \cite{AD} on pile sizes in patience sorting.  

{\bf Remark}: An analogue of part (1) of the corollary for $d$ dimensional flat Minkowski space is given in \cite{BG}. Combined with the techniques of \cite{DZ} this can be extend to domains in arbitrary Lorentzian manifolds. We plan to pursue this line of reasoning in future work. 
\section{Finer asymptotics for the uniform distribution}  
In this section we consider finer asymptotics for $ST(\bar{R})$ when the distribution $p$ 
on $C$ is uniform. 

Let $X_n$ be the random variable which counts the length of the 
longest increasing subsequence among $n$ uniformly chosen points 
in the unit square. Let $Y_n$ be the corresponding ``Poisonnized'' 
random variable which counts the length of the longest increasing 
subsequence for a Poisson process of intensity $n$ in the unit square. 
We recall that an intensity $n$ Poisson process will generate $k$ 
uniformly distributed points in the unit square with probability 
$e-n\frac{n^k}{k!}$ and that the restriction of the points to 
disjoint subsets are independent processes. The statistics of $Y_n$ have been studied in great detail, see \cite{BDJ}. 
The following result of M. Lowe and F. Merkl, see \cite{LM} theorem 1.2, computes the moderate deviations of $Y_n$ and $X_n$.
\begin{theo} \label{merkl}
Let $0<t_n<n^{\delta}$, for some $\delta <1/3$ be a sequence such that $t_n\longrightarrow \infty $. Let $l_n=t_nn^{1/6}+2\sqrt{n}$, then
\[ log(Pr(Y_n\geq l_n))\approx -(4/3)t_n^{3/2}\]
The same holds for $X_n$.
\end{theo}

We also need the following result of K. Johansson.
\begin{theo} 
\label{johansson}
(\cite{Jo}), Let $\delta >0$. Let $Z_n$ be the length of the longest increasing subsequence 
of an intensity $n$ Poisson process restricted to 
the diagonal stripe given by $\vert y-x\vert <n^{\delta-1/6}$ then 
$Pr(Y_n=Z_n)=1-e^{-n^{c(\delta )}}$, where $c(\delta )>0$. 
\end{theo}

Using these two results we prove the following finer estimate on $L_n$.
\begin{theo}
If $f(\theta )=c\theta $ and $p$ is the uniform distribution then, for all $\varepsilon >0$, w.h.p 
\[ Aln^{2/3}(n)<ST(\bar{R})-\sqrt{\frac{2}{c}}\sqrt{n}<Bln^{2/3}(n)\] with
$A=\frac{1}{4}(2c)^{-1/6}-\varepsilon $ and 
$B=3^{2/3}A+\varepsilon $.
\end{theo}
{\bf Proof}: For simplicity we assume that $c=1$. Consider the unit square $I^2\subset U$ with coordinates $(\theta , r)$
and the partial order relation $\leq_{ver}$ restricted from $U$. Let $\varepsilon >0$. Consider the 
set $I_{\varepsilon }$ given by $\varepsilon < \theta < 1-\varepsilon $ as a subset of either $C$ or $I^2$.
When $p$ is the uniform distribution or more generally of the form $p=p(r)$ we know by the proof of theorem \ref{estimate} that w.h.p the
maximal chain in $C$ will lie in a $\varepsilon $ neighborhood of a vertical line. 
Choose uniformly one of the lines say 
$\psi (r)=\theta_0$ for some $\theta_0$. Since no value of $\theta_0$ 
is more likely than any other value, we have with 
probability $1-4\varepsilon$ that $2\varepsilon <\theta_0 <1-2\varepsilon $ 
and thus the maximal chain is contained 
in $I_{\varepsilon}$. By definition any chain in $I^2$ w.r.t $\leq_{ver}$ is a chain in $C$ w.r.t $\leq_{C,ver}$.
When $\varepsilon <1/4$ then the relations $\leq_{C,ver}$ and $\leq_{U,ver}$ 
coincide on any $\varepsilon $ neighborhood 
of a vertical line. to see this note that if 
$P\leq_{C,ver}Q$ then by definition there is a $Q'$ such that $T(Q)=T(Q')$
and $r(Q')-r(P)\geq c(\vert t(Q')-t(p)\vert )$. 
By assumption $\vert t(Q)-t(P)\vert \leq 2\varepsilon <1/2$ hence if $Q'\neq Q$  
$\vert t(Q')-t(P)\vert\geq \vert t(Q')-t(Q)\vert -\vert t(Q)-t(P)\vert \geq 1-1/2
=1/2\geq \vert t(Q)-t(P)\vert $. We 
conclude that $P\leq_{U,ver}  Q$ as desired. Taking the limit as 
$\varepsilon $ goes to $0$ We see that w.h.p 
$M_{C,ver}(\bar{R})=M_{I^2,ver}(\bar{R})$. As a consequence of theorem \ref{combinatorial} we are reduced to finding finer estimates for
$M_{I^2,ver}(\bar{R})$, where the requests are uniformly chosen. Rotating $I^2$ by $\pi /4$ and then shifting and scaling
 the image 
we see that this is the same as finding $M_{K,inc}(\bar{R})$ where $K$ is a diamond with vertices 
$(0,-1/\sqrt{2}),(1/\sqrt{2},0),(-1/\sqrt{2},0),(0,1/\sqrt{2})$ and $\bar{R}$ are uniformly chosen points in $K$. Let $\delta_n \longrightarrow 0$
be positive. Let 
$x_i=\frac{i}{n^{1/6-\delta}\sqrt{2}}$ and consider the points 
$A_i=(-x_i,-1/\sqrt{2}+x_i)$,  $i=1,...,n^{1/6-\delta_n }$, which all lie on the bottom left boundary of $K$.
Let $B_i=A_i+(1/\sqrt{2},1/\sqrt{2})$, all of whom lie on the top right boundary. We define $S_i$ to be the square with edges parallel to the axis and opposing vertices $A_i$ and $B_i$. Let $w=\frac{1}{\sqrt{2}}n^{-1/6+\delta_n}$ and define 
$C_i=A_i+(0,w)$, 
$D_i=A_i+(w,0)$, $E_i=B_i-(0,w)$
and $F_i=B_i-(w,0)$. Let $H_i$ be the hexagon whose vertices are 
$A_i,B_i,C_i,D_i,E_i,F_i$. $H_i$ is precisely the $n^{-1/6+\delta}$
neighborhood of the diagonal in the square $S_i$. $H_i$ meets $H_{i+1}$ 
along the mutual
edge whose endpoints are $C_i=D_{i+1}$ and $F_i=E_{i+1}$. the interiors of the $H_i$ are disjoint.
The square $S_i$ has area $1/2$. By scaling, an intensity $n/2$ Poisson process on $S_i$ is equivalent 
to an intensity $n$ process
on the unit square. By theorem \ref{johansson} w.h.p the longest increasing subsequence 
in $H_i$ coincides with that of $S_i$ for
all $i$. We see that the longest increasing subsequence in $K$ 
dominates the random variable $\tilde{Y_n}=Max(Y_{n/2,i})$, $i=1,...,n^{1/6-\delta_n } $, where the $Y_{n/2,i}$ are i.i.d with 
distribution $Y_{n/2}$. Let $\varepsilon >0$.
Since the $Y_{n/2,i}$ are independent $\tilde{Y}>2\sqrt{n/2}+t_n(n/2)^{1/6}$ w.h.p if $Pr(Y_{n/2}>\sqrt{2n}+t_n(n/2)^{1/6})>n^{-1/6+\varepsilon }$.
by theorem \ref{merkl} the last inequality holds if $t_n=\frac{1}{4}ln^{2/3}(n)(1-\varepsilon ')$ for all $\varepsilon '>0$. Letting
$\varepsilon $ go to zero yields the lower bound.

To obtain the upper bound we let $z_i=\frac{i}{n^{1/3+\delta}\sqrt{2}}$, $i=1,...,n^{1/3+\delta }$. 
We let $P_i=(-z_i,-1/\sqrt{2}+z_i)$
and let $Q_i=P_i+(1/\sqrt{2},1/\sqrt{2})$
Let $u=(1/\sqrt{2})n^{-1/3+\delta}$, 
$T_i=P_i-(0,u)$ 
and $U_i=Q_i+(u,0)$. 
Let $R_{i,j}$ be the rectangle with sides parallel to the axis and whose opposite vertices are $T_i$ and $U_j$. It is easy to see
that $K$ is contained in $L=\cup_{i,j} R_{i,j}$ (in fact in the union of $R_{i,i}$). The $T_i$ form the set of minimal points
of the poset $(L, inc)$, likewise $U_j$ form the set of maximal points, consequently, every chain in $L$ is contained in $R_{i,j}$
for some pair $i,j$.
The area of $R_{i,i}$ is $vol(R_{i,j})=(1/\sqrt{2}+n^{-1/3-\delta })^2=\frac{1}{2}(1+2n^{-1/3-\delta}+n^{-2/3-2\delta})$. consider a pair
$i,j$ and let $k=\vert i-j\vert $. the area of $R_{i,j}$ is 
$vol(R_{i,j})=\frac{1}{2}(1+(k+1)n^{-1/3-\delta })(1-(k-1)n^{-1/3-\delta })=\frac{1}{2}(1+2n^{-1/3-\delta }-(k^2-1)n^{-2/3-2\delta })$. Let $\eta $ be any real satisfying $\delta <\eta $. We see that if $k>2n^{1/6+\eta }$ we have $vol(R_{i,j})<\frac{1}{2}(1 -n^{-1/3+2(\eta-\delta )})$. Applying theorem \ref{merkl} to $\frac{1}{2}(1-n^{-1/3+2(\eta -\delta )})$ we see that if $k>2n^{1/6+\eta }$ the probability that the 
longest increasing subsequence in 
$R_{i,j}$ has length greater than $\sqrt{2n}$ is at most $e^{-n^{\nu }}$ for some $\nu >0$. from our lower bound we conclude that
with high probability the longest increasing subsequence in $L$ is contained in some $R_{i,j}$ with $\vert i-j\vert <2n^{1/6+\eta }$.
The number of pairs satisfying this relation is at most $4n^{1/2+\delta +\eta }$
We also note that $vol(R_{i,j})<vol(R_{i,i})=\frac{1}{2}(1+2n^{-1/3-\delta}+n^{-2/3-2\delta })$ and hence the distribution of the length of the longest increasing subsequence in $R_{i,j}$ is dominated by that of $R_{i,i}$. Let $m=vol(R_{i,i})n$. We conclude that the 
distribution of the length of the longest increasing sequence in $L$ is dominated by $\bar{Y}=Max(Y_{m,i})$, $i=1,..,4n^{1/2+\delta +\eta }$. Using the union bound, applying theorem \ref{merkl} as before and letting $\delta ,\eta $ tend to zero we see that w.h.p $Pr(\bar{Y}<\sqrt{2n}+t_n(n/2)^{1/6})$
for $t_n=3^{2/3}\frac{1}{4}ln^{2/3}(n)(1-\varepsilon )$. Finally, our estimates can be depoissonized as in \cite{LM} or \cite{BDJ} section 6. 
{\em q.e.d}

{\it Acknowledgments}:
I would like to thank D. Berend, L. Sapir, S. Skiena and especially O. Zeitouni for
very helpful conversations.

{\bf Appendix: A very brief introduction to Lorentzian geometry}

We briefly consider some very basic material on Lorentzian geometry which is relevant
to our application. Since Lorentzian geometry is intimately tied to relativity theory and
in fact was invented to model it, we also give a brief account of the physical interpretation
of various geometrical notions in the theory. For a concise but much more comprehensive mathematical treatment the reader is referred to
\cite{Pe} and for more on the physics to \cite{Ca}.

A Lorentzian metric on a domain $D\subset \bf{R}^n$ with local coordinates $x_1,...,x_n$ 
is given by a $C^{\infty }$ mapping from $D$ to 2-forms 
\begin{equation}
ds^2=\sum_i\sum_jg_{ij}dx_idx_j
\end{equation}
Were $g_{ij}$ is a symmetric non singular $n$ by $n$ matrix with a single positive eigenvalue. 
The form acts like the square of a ``distance'' function
between nearby points. 
A space-time is a manifold (with boundary) with a global Lorentzian metric.

{\bf Example}: Consider $\bf{R}^{n+1}$, denote the first coordinate which we index 
as the 0'th coordinate 
by $t$ and the other coordinates by $x_1,...,x_n$.
We attach to all points in $\bf{R}^{n+1}$ the fixed diagonal matrix $g$ with entries $(1,-1,-1,...,-1)$. 
We have $ds^2=dt^2-dx_1^2-...-dx_n^2$.
$\bf{R}^{n+1}$ equipped with this constant Lorentzian metric is 
known as flat Minkowski space. It was introduced by Minkowski 
to provide a mathematical (geometrical) model for special relativity theory and 
in particular the notion of space-time as opposed to 
space and time separately.  

Given the metric $g$, we say that a tangent vector $v$ is {\it time-like} if $vgv^t>0$, 
{\it null} if $vgv^t=0$ and {\it space-like} if $vgv^t<0$. We say that $v$ and $w$ are orthogonal
if $vgw^t=0$. The metric $g$ also induces a volume form which in local coordinates is given by the formula
\begin{equation}
\label{volume}
vol=\sqrt{\vert det(g)\vert }dx_1...dx_{n+1}
\end{equation}
When a curve $\phi (u)$ is piecewise differentiable the length of $\phi $ 
which we denote by $\ell (\phi )$ is given by the formula $\int_a^b \sqrt{\phi '(u)g\phi '(u)^t}du$
where $\phi '(u)=(x_1'(u),...,x_n'(u))$ is the tangent of $\phi $ at $u$. 
$\ell (\phi )$ can only be defined if $\phi '(u)g\phi '(u)^t\geq 0$ for all $u$.

We see that length can only 
be defined for curves whose tangents are either time-like or light-like. 
The set of time-like vectors has two cone shaped components. 
When there are no closed time-like curves we can define a notion of past and future
, causality, via a consistent choice at each point of one of the components, 
that is, a continuous vector field of non vanishing time-like vectors. The vectors in the chosen 
component will be called future pointing. We can then define a partial
order on the points of the space-time by declaring that $A\leq_gB$ if there is a curve beginning in $A$ and 
ending in $B$ with future pointing (time-like) tangents. Physically, the partial order corresponds to the notion of causality, past-future relation. We note that scaling the metric by a (point dependent)
scalar function does not affect this partial order. Such changes to the metric are called
{\it conformal}.

Relativity theory models gravity using a Lorentzian metric on a manifold (the universe). The motion of a particle (small body) 
with positive mass is described by a curve with future pointing time-like tangents.
Such a curve is known as a {\it world line}. The 
length of the curve is the time which passes when measured by a clock attached 
to the particle, this is known in 
physics as {\it proper time}. 
Massless particles travel along curves
with null tangents. When no force other than gravity is exerted on the particle (free falling particle), 
the trajectory of the particle can be divided into small segments so that each segment 
maximizes proper time (curve length) among all paths with the same endpoints. 
Curves with this property are known as {\it geodesics}. Geodesics are always differentiable curves.
In flat Minkowski space which models a world without gravity, such curves are straight lines. 
This is a restatement 
of Galileo's principle (the first law of mechanics) that a 
free particle will travel in a straight line at constant speed.

In two dimensional space-time, if $g$ is a Lorentzian metric than so is $-g$. The causal structures induced by $g$ and $-g$
are complementary in the sense that $v$ is space-like w.r.t $g$ iff it is time-like w.r.t $-g$. Consequently, we can use $-g$ to measure
the length of curves which are space-like w.r.t $g$.
\end{document}